\documentclass{amsart}

\usepackage{pstricks,pst-node,graphicx,amsmath}

\newtheorem{theorem}{Theorem}[section]

\newtheorem{proposition}[theorem]{Proposition}
\newtheorem{corollary}[theorem]{Corollary}

\theoremstyle{definition}
\newtheorem{definition}[theorem]{Definition}

\theoremstyle{remark}
\newtheorem{remark}[theorem]{Remark}

\numberwithin{equation}{section}

\newcommand{\SYT}{\ensuremath\mathrm{SYT}}

\newcommand{\GL}{\ensuremath\mathrm{GL}}
\newcommand{\SL}{\ensuremath\mathrm{SL}}
\newcommand{\Sn}{\ensuremath\mathcal{S}_n}

\newcommand{\wf}{\ensuremath\tilde{f}}
\newcommand{\we}{\ensuremath\tilde{e}}

\newcommand{\ipo}{\ensuremath i\!+\!1}
\newcommand{\imo}{\ensuremath i\!-\!1}
\newcommand{\imt}{\ensuremath i\!-\!2}
\newcommand{\imh}{\ensuremath i\!-\!3}

\newcommand{\nmo}{\ensuremath n\!-\!1}

\newcommand{\G}{\ensuremath\mathcal{G}}
\newcommand{\X}{\ensuremath\mathcal{X}}

\newcommand{\ial}{\ensuremath\stackrel{i}{\longleftarrow}}
\newcommand{\iar}{\ensuremath\stackrel{i}{\longrightarrow}}

\newcommand{\reffig}[1]{Figure \ref{fig:#1}}
\newcommand{\refeq}[1]{equation (\ref{eqn:#1})}

\newcommand{\refthm}[1]{Theorem \ref{thm:#1}}

\newcommand{\refcor}[1]{Corollary \ref{cor:#1}}

\newcommand{\refsec}[1]{Section \ref{sec:#1}}

\newlength{\hsp}
\newlength{\vsp}
\newlength{\vspi}
\newcommand{\cci}  {c@{\hskip  \hsp}} 
\newcommand{\ccii} {c@{\hskip 2\hsp}} 
\newcommand{\cciii}{c@{\hskip 3\hsp}} 

\newcommand{\rn}{\rnode}

\newcommand{\B}{\bullet}

\newlength\cellsize 

\savebox2{%
\begin{picture}(10,10)
\put(0,0){\line(1,0){10}}
\put(0,0){\line(0,1){10}}
\put(10,0){\line(0,1){10}}
\put(0,10){\line(1,0){10}}
\end{picture}}

\newcommand\cellify[1]{\def\thearg{#1}\def\nothing{}%
\ifx\thearg\nothing
\vrule width0pt height\cellsize depth0pt\else
\hbox to 0pt{\usebox2\hss}\fi%
\vbox to 10\unitlength{
\vss
\hbox to 10\unitlength{\hss$#1$\hss}
\vss}}

\newcommand\tableau[1]{\vtop{\let\\=\cr
\setlength\baselineskip{-10000pt}
\setlength\lineskiplimit{10000pt}
\setlength\lineskip{0pt}
\halign{&\cellify{##}\cr#1\crcr}}}

\newcommand{\stab}[3]{\begin{array}{c}\rnode{#1}{\tableau{#2}}\\\rnode{#1#1}{_{#3}}\end{array}}

\begin{document}


\title[Combinatorial Schur-Weyl duality]{A combinatorial realization
  of Schur-Weyl duality via crystal graphs and dual equivalence graphs}

\author[S. Assaf]{Sami H. Assaf}
\address{Department of Mathematics, University of Pennsylvania,
Philadelphia, PA 19104}
\email{sassaf@math.mit.edu}
\thanks{Work supported in part by an NSF Mathematical Sciences
  Postdoctoral Research Fellowship DMS-0703567.}

\subjclass[2000]{%
Primary   05e10; 
Secondary 05e15, 
          20c30, 
}



\begin{abstract}
  For any polynomial representation of the special linear group, the
  nodes of the corresponding crystal may be indexed by semi-standard
  Young tableaux. Under certain conditions, the standard Young
  tableaux occur, and do so with weight $0$. Standard Young tableaux
  also parametrize the vertices of dual equivalence graphs. Motivated
  by the underlying representation theory, in this paper, we explain
  this connection by giving a combinatorial manifestation of
  Schur-Weyl duality. In particular, we put a dual equivalence graph
  structure on the $0$-weight space of certain crystal graphs,
  producing edges combinatorially from the crystal edges. The
  construction can be expressed in terms of the local
  characterizations given by Stembridge for crystal graphs and the
  author for dual equivalence graphs.
\end{abstract}

\maketitle

%
\section{Introduction}
%
\label{sec:introduction}

Schur-Weyl duality \cite{Weyl1939} is a powerful tool in the study of
irreducible representations of the classical groups. The set-up is as
follows. Let $G$ be a complex, reductive Lie group, and let $V$ be an
irreducible finite-dimensional complex representation of $G$. Let
$V^0$ be the space of vectors of weight $0$, i.e. vectors which are
fixed by a maximal torus. Then the Weyl group $W$ acts naturally on
$V^0$. The case when $G = \SL_n$, and so then the Weyl group is
isomorphic to $\Sn$, was studied by Gutkin \cite{Gutkin1973} who
showed that if $V$ is an irreducible component of the $n$th tensor
power of the standard $n$-dimensional representation of $\SL_n$, then
the $\Sn$-module $V^0$ is irreducible and all irreducible
$\Sn$-modules occur in this way. This was also observed by Kostant
\cite{Kostant1976}. Though there are descriptions of the Weyl group
action on $0$-weight spaces for $\SL_n$, there is no general rule for
how to decompose the action into irreducible representations of $\Sn$.
This paper proposes a new approach to this problem by
combinatorializing the representations using crystal graphs
\cite{Kashiwara1990,Kashiwara1991} for representations of $\SL_n$ and
dual equivalence graphs \cite{Assaf2007-2} for representations of
$\Sn$.

Kashiwara introduced in \cite{Kashiwara1990,Kashiwara1991} the notion
of crystal bases in his study of the representation theory of
quantized universal enveloping algebras at $q=0$. The theory of
canonical bases, developed independently by Lusztig
\cite{Lusztig1990,Lusztig1990-2}, studies the same problem from a
different viewpoint, though for the purposes of this paper, we focus
solely on the crystal approach. A crystal graph is a directed, colored
graph with vertex set given by the crystal basis and directed edges
given by deformations of the Chevalley generators $e_i$ and $f_i$. The
combinatorial structure of crystal graphs encodes important
information for studying the representations; for instance, knowing
the crystal immediately gives a formula for the character as well as
tensor product and branching rules for the corresponding
representation. For the quantum group $U_{q}(\mathfrak{sl}_n)$, the
crystal basis can be indexed by semi-standard Young tableaux, and
there is an explicit combinatorial construction of the crystal graph
on tableaux \cite{KaNa1994,Littelmann1995}.

Dual equivalence graphs were first introduced in \cite{Assaf2007-2} as
a combinatorial tool for studying functions expressed in terms of
quasi-symmetric functions. In particular, in \cite{Assaf2007-2} they
are used to give a purely combinatorial proof of symmetry and Schur
positivity of Lascoux-Leclerc-Thibon polynomials and Macdonald
polynomials. The original motivation for these graphs, based on ideas
of Haiman \cite{Haiman2005}, was to mimic existing crystal-theoretic
proofs used to establish positivity results for functions expressed in
terms of monomials. Given this, a natural problem which this paper
addresses is to explore the connections between crystal graphs and
dual equivalence graphs. Combinatorially, the vertices of a dual
equivalence graph are standard Young tableaux of a given shape, which
also index a basis for irreducible representations of $\Sn$. When
$\lambda$ is a partition of $n$, these also index the $0$-weight nodes
in the crystal graph of the irreducible $\SL_n$ module with highest
weight $\lambda$. From this observation we develop a correspondence
between crystals of irreducible representations of $\SL_n$ of degree
$n$ and dual equivalence graphs.

This paper is organized as follows. \refsec{crystals} contains an
exposition of type A crystal graphs. We begin by defining the
combinatorial structure on semi-standard Young tableaux, and then
present Stembridge's local characterization for crystals
\cite{Stembridge2003} which will be essential for establishing our
main result. \refsec{degs} gives the analogous exposition for dual
equivalence graphs, beginning with the combinatorial description in
terms of standard Young tableaux before presenting the local
characterization from \cite{Assaf2007-2}. The correspondence is laid
out in \refsec{correspondence}, where we give two separate proofs of
the main result using the tableaux descriptions (\refthm{tableaux-pf})
and the local characterizations (\refthm{axiom-pf}). Finally, in
\refsec{extensions}, we lay out how these methods may be extended to
any irreducible representation of $\SL_n$ with a nontrivial $0$-weight
space, and we also suggest how this approach may be used to define
dual equivalence graphs for other types.

\begin{center}
{\sc Acknowledgements}
\end{center}

The author is grateful to Mark Haiman for the idea to construct the
standard dual equivalence graph $\G_{\lambda}$ and for suggesting
investigating the connections with crystal graphs. The author also
thanks Ian Grojnowski, Mark Haiman and Monica Vazirani for helpful
conversations about crystal graphs and $0$-weight spaces.

%
\section{Crystal graphs}
%
\label{sec:crystals}

In order to construct the crystal graph of an irreducible
representation of $\SL_n$, we must define Kashiwara's crystal
operators in this setting. We depart somewhat from the original
descriptions developed independently by Kashiwara and Nakashima
\cite{KaNa1994} and Littelmann \cite{Littelmann1995} in favor of the
presentation of the Littlewood-Richardson rule in
\cite{Macdonald1995}.

For a word $w$ of length $n$, i.e. $w \in \mathbb{N}^{n}$, define
$m_i(w,r)$ to be the number of $i$'s occurring in $w_{r} w_{r+1}
\cdots w_{n}$ minus the number of $\imo$'s occurring in $w_{r} w_{r+1}
\cdots w_{n}$. Then let $m_i(w)$ be the maximum $m_i(w,r)$ which
occurs for any $r = 1, \ldots, n$.  Observe that if $m_i(w) > 0$ and
$w_p$ is the rightmost occurrence of this maximum, i.e. $m_i(w,p) =
m_i(w)$ and $m_i(w,q) < m_i(w)$ for $q > p$, then $w_p = i$. The
result of applying the crystal operator $\we_i$ to a semi-standard
tableau $T$ may be described as follows. Let $w_T$ be the reading word
for $T$. Then if $m_i(w_T) \leq 0$, $\we_i(T)=0$; otherwise $\we_i(T)$
is the result of changing the $i$ which is the rightmost occurrence of
$m_i(w_T)$ to an $\imo$. As $\we_i$ is invertible when nonzero, define
$\wf_i$ by requiring $\we_i$ and $\wf_i$ to be inverses of one another
whenever the image is not $0$.

\begin{center}
  \begin{figure}[ht]
    \begin{displaymath}
      \begin{array}{\ccii\ccii\ccii\ccii \ccii \ccii\ccii\ccii\ccii}
        & & & & \rn{e1}{\tableau{3 & 3 \\ 2 & 2}} & & & & \\[1.5\vsp]
        & & \rn{c2}{\tableau{2 & 3 \\ 1 & 2}} &
        \rn{d2}{\tableau{3 & 3 \\ 1 & 2}} & &
        \rn{f2}{\tableau{3 & 4 \\ 2 & 2}} &
        \rn{g2}{\tableau{3 & 4 \\ 2 & 3}} & & \\[1.5\vsp]
        \rn{a3}{\tableau{2 & 2 \\ 1 & 1}} &
        \rn{b3}{\tableau{2 & 3 \\ 1 & 1}} &
        \rn{c3}{\tableau{3 & 3 \\ 1 & 1}} &
        \rn{d3}{\tableau{2 & 4 \\ 1 & 2}} &
        \begin{array}{c}
          \rn{e3l}{\tableau{2 & 4 \\ 1 & 3}} \\[\vsp]
          \rn{e3r}{\tableau{3 & 4 \\ 1 & 2}}
        \end{array} &
        \rn{f3}{\tableau{3 & 4 \\ 1 & 3}} &
        \rn{g3}{\tableau{4 & 4 \\ 2 & 2}} &
        \rn{h3}{\tableau{4 & 4 \\ 2 & 3}} & 
        \rn{i3}{\tableau{4 & 4 \\ 3 & 3}} \\[2.5\vsp]
        & & \rn{c4}{\tableau{2 & 4 \\ 1 & 1}} &
        \rn{d4}{\tableau{3 & 4 \\ 1 & 1}} & &
        \rn{f4}{\tableau{4 & 4 \\ 1 & 2}} &
        \rn{g4}{\tableau{4 & 4 \\ 1 & 3}} & & \\[1.5\vsp]
        & & & & \rn{e5}{\tableau{4 & 4 \\ 1 & 1}} & & & &
      \end{array}
      \psset{linewidth=.1ex,nodesep=3pt}
      \everypsbox{\scriptstyle}
      \ncline[linecolor=blue,linestyle=dashed]{->} {b3}{c2} 
      \ncline[linecolor=blue,linestyle=dashed]{->} {c3}{d2} 
      \ncline[linecolor=blue,linestyle=dashed]{->} {c4}{d3} 
      \ncline[linecolor=blue,linestyle=dashed]{->} {d2}{e1} 
      \ncline[linecolor=blue,linestyle=dashed]{->} {d4}{e3r}
      \ncline[linecolor=blue,linestyle=dashed]{->} {e3r}{f2}
      \ncline[linecolor=blue,linestyle=dashed]{->} {e5}{f4} 
      \ncline[linecolor=blue,linestyle=dashed]{->} {f3}{g2} 
      \ncline[linecolor=blue,linestyle=dashed]{->} {f4}{g3} 
      \ncline[linecolor=blue,linestyle=dashed]{->} {g4}{h3} 
      \ncline[linecolor=red]{->}  {a3}{b3} 
      \ncline[linecolor=red]{->}  {b3}{c3} 
      \ncline[linecolor=red]{->}  {c2}{d2} 
      \ncline[linecolor=red]{->}  {c4}{d4} 
      \ncline[linecolor=red]{->}  {d3}{e3l}
      \ncline[linecolor=red]{->}  {e3l}{f3}
      \ncline[linecolor=red]{->}  {f2}{g2} 
      \ncline[linecolor=red]{->}  {f4}{g4} 
      \ncline[linecolor=red]{->}  {g3}{h3} 
      \ncline[linecolor=red]{->}  {h3}{i3} 
      \ncline[linecolor=green,linewidth=.3ex]{->}{b3}{c4} 
      \ncline[linecolor=green,linewidth=.3ex]{->}{c2}{d3} 
      \ncline[linecolor=green,linewidth=.3ex]{->}{c3}{d4} 
      \ncline[linecolor=green,linewidth=.3ex]{->}{d2}{e3r}
      \ncline[linecolor=green,linewidth=.3ex]{->}{d4}{e5} 
      \ncline[linecolor=green,linewidth=.3ex]{->}{e1}{f2} 
      \ncline[linecolor=green,linewidth=.3ex]{->}{e3r}{f4}
      \ncline[linecolor=green,linewidth=.3ex]{->}{f2}{g3} 
      \ncline[linecolor=green,linewidth=.3ex]{->}{f3}{g4} 
      \ncline[linecolor=green,linewidth=.3ex]{->}{g2}{h3} 
    \end{displaymath}
    \caption{\label{fig:X22}The crystal graph $\X_{(2,2)}^{4}$, with
      edges $\wf_1{\blue \nearrow}$, $\wf_2{\red
        \rightarrow}$, $\wf_3{\green \searrow}$.}
  \end{figure}
\end{center}

For a partition $\lambda$ with at most $n$ rows, the crystal graph
$\X_{\lambda}^{n}$ is the directed, colored graph on semi-standard
tableaux of shape $\lambda$ with entries from $[n]$ defined as
follows. For $i=1,\ldots,\nmo$, if $m_i(w_T) > 0$, define a directed
$i$-edge from $\we_i(T)$ to $T$. For example, \reffig{X22} gives the
crystal graph $\X_{(2,2)}^{4}$.

In \cite{Stembridge2003}, Stembridge gives a local characterization of
crystal graphs which arise from representations for simply-laced
types. In order to define Stembridge's axioms, we must first introduce
some notation associated with a directed, colored graph. Say that $\X$
has {\em degree $n$} if the largest color for an edge is $\nmo$. 

We have simplified things from the presentation in
\cite{Stembridge2003} since we are considering only the type A
case. The notations below may not be well-defined in general, but will
be in the contexts in which they are used, i.e. for graphs satisfying
axioms P1 and P2 below.

If $x {\blue {\blue \ial}} y$ (resp. $x {\blue \iar} z$), write $y=E_i
x$ (resp. $z = F_i x$). The {\em $i$-string through $x$} is the
maximal path
$$
F_i^{-d} x {\blue \iar} \cdots {\blue \iar} F_i^{-1} x {\blue \iar} x
{\blue \iar} F_i x {\blue \iar} \cdots {\blue \iar} F_i^{r} x .
$$
In this case we write $\delta(x,i) = -d$ and $\varepsilon(x,i) =
r$. Finally, we have the following differences whenever $E_i, F_i$ is
defined at $x$.
\begin{displaymath}
  \begin{array}{rclcrcl}
    \Delta_i \delta(x,j) & = & \delta(E_i x, j) - \delta(x,j), & &
    \nabla_i \delta(x,j) & = & \delta(x, j) - \delta(F_i x,j), \\
    \Delta_i \varepsilon(x,j) & = & \varepsilon(E_i x,j) - \varepsilon(x,j), & &
    \nabla_i \varepsilon(x,j) & = & \varepsilon(x,j) - \varepsilon(F_i x,j).
  \end{array}
\end{displaymath}

\begin{definition}[\cite{Stembridge2003}]
  A directed, colored graph $\X$ is {\em regular} if the following
  hold:
  \begin{itemize}
  \item[(P1)] all monochromatic directed paths have finite length;

  \item[(P2)] for every vertex $x$, there is at most one edge $x {\blue \ial}
    y$ and at most one edge $x {\blue \iar} z$;

  \item[(P3)] assuming $E_i x$ is defined, $\Delta_i \delta(x,j) +
    \Delta_i \varepsilon(x,j) = \left\{\begin{array}{rl}
        2 & \;\mbox{if}\;\; j=i \\
        -1 & \;\mbox{if}\;\; j=i\pm 1 \\
        0 & \;\mbox{if}\;\; |i-j|\geq 2
      \end{array}\right.$;
  \item[(P4)] assuming $E_i x$ is defined, $\Delta_i \delta(x,j),
    \Delta_i \varepsilon(x,j) \leq 0$ for $j \neq i$;

  \item[(P5)] $\Delta_i \delta(x,j) = 0$ $\Rightarrow$ $E_iE_j x = E_jE_i
    x = y$ and $\nabla_j \varepsilon(y, i) = 0$; \\
    $\nabla_i \varepsilon(x,j) = 0$ $\Rightarrow$ $F_iF_j x = F_jF_i
    x = y$ and $\Delta_j \delta(y, i) = 0$;

  \item[(P6)] $\Delta_i \delta(x,j) = \Delta_j \delta(x,i) = -1$
    $\Rightarrow$ $E_iE_{j}^{2}E_i x = E_jE_{i}^{2}E_j x = y$ and
    $\nabla_i \varepsilon(y,j) = \nabla_j \varepsilon(y,i)=-1$; \\
    $\nabla_i \varepsilon(x,j) = \nabla_j \varepsilon(x,i) = -1$
    $\Rightarrow$ $F_iF_{j}^{2}F_i x = F_jF_{i}^{2}F_j x = y$ and
    $\Delta_i \delta(y,j) = \Delta_j \delta(y,i)=-1$.
  \end{itemize}
  \label{defn:regular}
\end{definition}

\begin{figure}[ht]
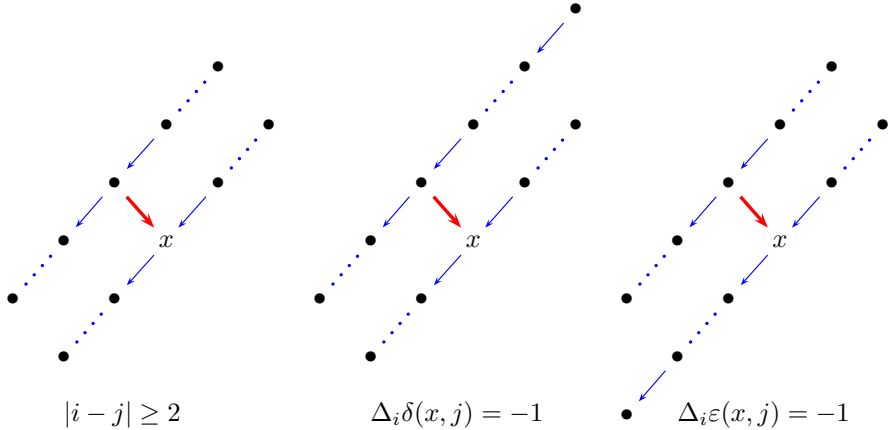

  \begin{center}
    \begin{displaymath}
      \begin{array}{\cci\cci\cci\cci\cci\cci 
          \cci\cci\cci\cci\cci\cci \cci\cci\cci\cci\cci\cci}
        & & & & & & & & & & & \rn{a12}{\B} & & & & & & \\[.7\vsp]
        & & & &  \rn{b5}{\B} & &
        & & & & \rn{b11}{\B} & & 
        & & & & \rn{a17}{\B} & \\[.7\vsp]
        & & & \rn{c4}{\B}  & & \rn{c6}{\B}  &
        & & & \rn{c10}{\B} & & \rn{c12}{\B} &
        & & & \rn{b16}{\B} & & \rn{b18}{\B} \\[.7\vsp]
        & & \rn{d3}{\B}  & & \rn{d5}{\B}  & &
        & & \rn{d9}{\B}  & & \rn{d11}{\B} & &
        & & \rn{c15}{\B} & & \rn{c17}{\B} & \\[.7\vsp]
        & \rn{e2}{\B}  & & \rn{e4}{x}  & & &
        & \rn{e8}{\B}  & & \rn{e10}{x} & & &
        & \rn{d14}{\B} & & \rn{d16}{x} & & \\[.7\vsp]
        \rn{f1}{\B}  & & \rn{f3}{\B}  & & & &
        \rn{f7}{\B}  & & \rn{f9}{\B}  & & & &
        \rn{e13}{\B} & & \rn{e15}{\B} & & & \\[.7\vsp]
        & \rn{g2}{\B}  & & & & &
        & \rn{g8}{\B}  & & & & &        
        & \rn{f14}{\B} & & & & \\[.7\vsp]
        & \makebox[0pt][l]{$|i-j| \geq 2$} & & & & & 
        & \makebox[0pt][l]{$\Delta_i \delta(x,j) = -1$} & & & & &
        \rn{g13}{\B}  & \makebox[0pt][l]{$\Delta_i
          \varepsilon(x,j) = -1$} & & & &  
      \end{array}
      \psset{linewidth=.1ex,nodesep=3pt}
      \everypsbox{\scriptstyle}
      \ncline[linecolor=blue,linestyle=dotted,linewidth=.3ex] {b5}{c4}
      \ncline[linecolor=blue]{->}              {c4}{d3} 
      \ncline[linecolor=blue]{->}              {d3}{e2} 
      \ncline[linecolor=blue,linestyle=dotted,linewidth=.3ex] {e2}{f1}
      \ncline[linewidth=.3ex,linecolor=red]{->}  {d3}{e4}  
      \ncline[linecolor=blue,linestyle=dotted,linewidth=.3ex] {c6}{d5}
      \ncline[linecolor=blue]{->}              {d5}{e4} 
      \ncline[linecolor=blue]{->}              {e4}{f3} 
      \ncline[linecolor=blue,linestyle=dotted,linewidth=.3ex] {f3}{g2}
      \ncline[linecolor=blue]{->}              {a12}{b11} 
      \ncline[linecolor=blue,linestyle=dotted,linewidth=.3ex] {b11}{c10}
      \ncline[linecolor=blue]{->}              {c10}{d9}  
      \ncline[linecolor=blue]{->}              {d9} {e8}  
      \ncline[linecolor=blue,linestyle=dotted,linewidth=.3ex] {e8} {f7}
      \ncline[linewidth=.3ex,linecolor=red]{->}  {d9}{e10}  
      \ncline[linecolor=blue,linestyle=dotted,linewidth=.3ex] {c12}{d11}
      \ncline[linecolor=blue]{->}              {d11}{e10} 
      \ncline[linecolor=blue]{->}              {e10}{f9}  
      \ncline[linecolor=blue,linestyle=dotted,linewidth=.3ex] {f9} {g8}
      \ncline[linecolor=blue,linestyle=dotted,linewidth=.3ex] {a17}{b16}
      \ncline[linecolor=blue]{->}              {b16}{c15} 
      \ncline[linecolor=blue]{->}              {c15}{d14} 
      \ncline[linecolor=blue,linestyle=dotted,linewidth=.3ex] {d14}{e13}
      \ncline[linewidth=.3ex,linecolor=red]{->}  {c15}{d16}  
      \ncline[linecolor=blue,linestyle=dotted,linewidth=.3ex] {b18}{c17}
      \ncline[linecolor=blue]{->}              {c17}{d16} 
      \ncline[linecolor=blue]{->}              {d16}{e15} 
      \ncline[linecolor=blue,linestyle=dotted,linewidth=.3ex] {e15}{f14}
      \ncline[linecolor=blue]{->}              {f14}{g13} 
    \end{displaymath}
    \caption{\label{fig:P3P4}An illustration of axioms P3 and P4,
      where $F_j {\blue \swarrow}$, $F_i {\red \searrow}$.}
  \end{center}
\end{figure}  

\begin{figure}[ht]
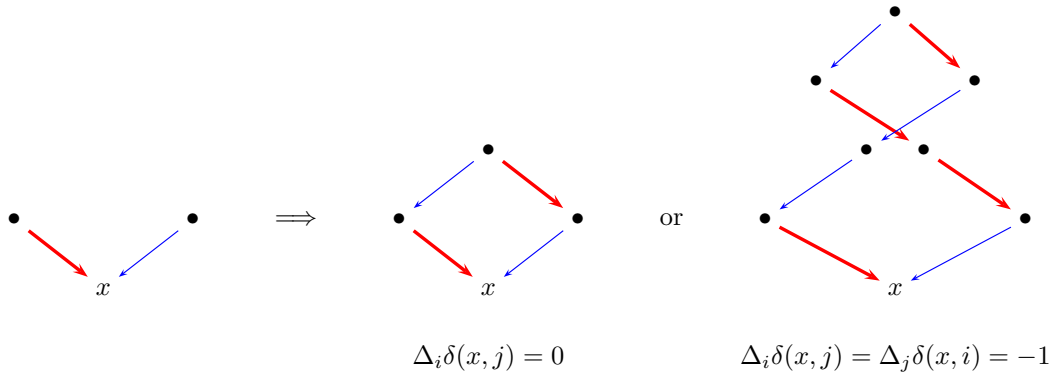

  \begin{center}
    \begin{displaymath}
      \begin{array}{\ccii\ccii\ccii \ccii \ccii\ccii\ccii \ccii
          \cci\cci\cci\cci\cci}  
        & & & & & & & & & & \rn{a11}{\B} & & \\[\vsp]
        & & & & & & & & & \rn{b10}{\B} & & \rn{b12}{\B} & \\[\vsp]
        & & & & & \rn{c6}{\B} & & & & & \rn{c11l}{\B} \;\;\;\;\;\;
        \rn{c11r}{\B} & & \\[\vsp] 
        \rn{d1}{\B} & & \rn{d3}{\B} & \rn{d4}{\Longrightarrow} &
        \rn{d5}{\B} & & \rn{d7}{\B} & \rn{d8}{\mbox{or}} & \rn{d9}{\B}
        & & & & \rn{d13}{\B} \\[\vsp]
        & \rn{e2}{x} & & & & \rn{e6}{x} & & & & & \rn{e11}{x} & & \\[\vsp]
        & & & & & \makebox[0pt][c]{$\Delta_i \delta(x,j) = 0$}
        & & & & & \makebox[0pt][c]{$\Delta_i \delta(x,j) = \Delta_j
          \delta(x,i) = -1$} & & 
      \end{array}
      \psset{linewidth=.1ex,nodesep=3pt}
      \everypsbox{\scriptstyle}
      \ncline[linecolor=red,linewidth=.3ex]{->} {d1}{e2}  
      \ncline[linecolor=blue]{->} {d3}{e2}  
      \ncline[linecolor=blue]{->} {c6}{d5}  
      \ncline[linecolor=red,linewidth=.3ex]{->} {c6}{d7}  
      \ncline[linecolor=red,linewidth=.3ex]{->} {d5}{e6}  
      \ncline[linecolor=blue]{->} {d7}{e6}  
      \ncline[linecolor=blue]{->} {a11}{b10}  
      \ncline[linecolor=red,linewidth=.3ex]{->} {a11}{b12}  
      \ncline[linecolor=red,linewidth=.3ex]{->} {b10}{c11r} 
      \ncline[linecolor=blue]{->} {b12}{c11l} 
      \ncline[linecolor=blue]{->} {c11l}{d9}  
      \ncline[linecolor=red,linewidth=.3ex]{->} {c11r}{d13} 
      \ncline[linecolor=red,linewidth=.3ex]{->} {d9}{e11}   
      \ncline[linecolor=blue]{->} {d13}{e11}  
    \end{displaymath}
    \caption{\label{fig:P5P6}An illustration of axioms P5 and P6,
      where $F_j {\blue \swarrow}$, $F_i {\red \searrow}$.}
  \end{center}
\end{figure}

Axioms P3 and P4 dictate how the maximal length of a $j$-string
differs between $x$ and $E_i x$; see \reffig{P3P4}. When $|i-j| \geq
2$, there is no change, but when $j = i\pm 1$, the length must change
by $1$. Axioms P5 and P6 give information about how edges with
different labels interact; see \reffig{P5P6}. When $|i-j| \geq 2$, the
conditions of P5 will always be satisfied, though when $j = i\pm 1$,
either P5 or P6 could hold.

The following theorem shows that regular graphs corresponds precisely
to the crystal graphs of representations in type A. This result allows
us to combinatorialize the problem of studying representations of
$\SL_n$ by instead studying regular graphs of degree $n$.

\begin{theorem}[\cite{Stembridge2003}]
  Each $\X_{\lambda}^{n}$ is a regular graph, and every connected
  component of a regular graph is isomorphic to $\X_{\lambda}^{n}$ for
  some $\lambda, n$.
\label{thm:structure-crystal}
\end{theorem}

The theorem is proved by showing that Littelmann's path operators
\cite{Littelmann1994} generate regular graphs, since the Path Model is
known to generate the $\X_{\lambda}^{n}$. As remarked in
\cite{Stembridge2003}, it would be nice to have a graph-theoretic
proof of \refthm{structure-crystal} which bypasses the Path Model.

%
\section{Dual equivalence graphs}
%
\label{sec:degs}

Since the characters of irreducible representations of $\SL_n$ are
Schur polynomials which form a basis for symmetric polynomials,
crystal graphs can be used to prove that certain polynomials expressed
in terms of monomials are symmetric and Schur positive. Motivated by
this approach, Haiman \cite{Haiman2005} suggested defining a graph
structure on standard tableaux using the dual equivalence relation as
a means of establishing the symmetry and Schur positivity of
polynomials expressed in terms of quasi-symmetric functions. The
result of this idea is the theory of dual equivalence graphs developed
in \cite{Assaf2007-2}.

Following Haiman's idea, we use the dual equivalence relation to
construct a graph whose vertices are given by standard tableaux and
whose connected components are indexed by partitions.

\begin{definition}[\cite{Haiman1992}]
  An {\em elementary dual equivalence} on three consecutive letters,
  say $\imo,i,\ipo$, of a permutation is given by switching the outer
  two letters whenever the middle letter is not $i$:
  \begin{center}
    $\displaystyle{%
      \cdots\; i \;\cdots\;i\pm 1\;\cdots\;i\mp 1\;\cdots \equiv^{*}
      \cdots\;i\mp 1\;\cdots\;i\pm 1\;\cdots\; i \;\cdots
    }$
  \end{center}
\label{defn:ede}
\end{definition}

The standard dual equivalence graph associated to $\lambda$, denoted
$\G_{\lambda}$, is the vertex-signed, edge-colored graph constructed
in the following way. The vertices are given by all standard Young
tableaux of shape $\lambda$. Two tableaux $T$ and $U$ are connected by
an $i$-colored edges whenever their reading words, $w_T$ and $w_U$,
differ by an elementary dual equivalence for $\imo,i,\ipo$.  To each
vertex $T$ we associate the signature $\sigma(T)$ which indicates the
descent set of the tableau:
\begin{equation}
  \sigma(T)_{i} \; = \; \left\{ 
    \begin{array}{ll}
      +1 & \; \mbox{if $i$ appears to the left of $\ipo$ in $w_T$} \\
      -1 & \; \mbox{if $\ipo$ appears to the left of $i$ in $w_T$}
    \end{array} \right. .
\label{eqn:sigma-T}
\end{equation}
As the original motivation for these graphs was to prove Schur
positivity, the role of the signatures is to index which
quasi-symmetric function should be associated to each vertex.

\begin{figure}[ht]
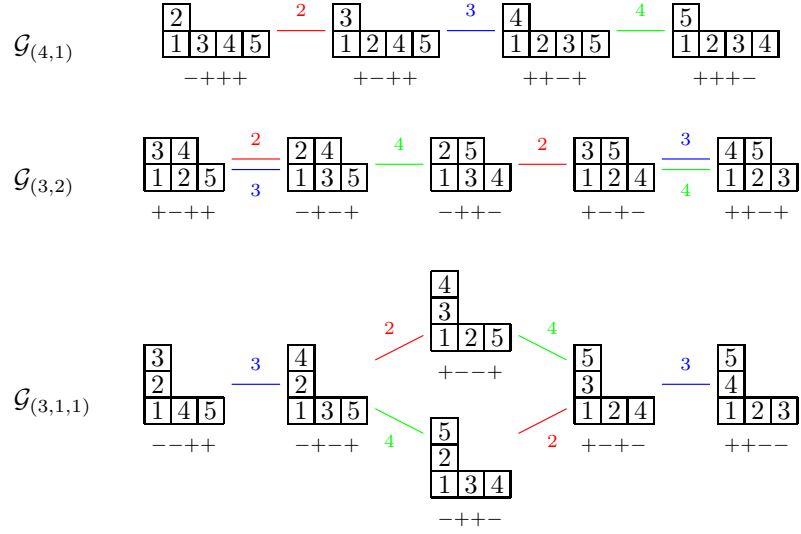

  \begin{center}
    \begin{displaymath}
      \begin{array}{lc}
        \G_{(4,1)} &
      \begin{array}{\cci \cci \cci \cci}
        \stab{h}{2 \\ 1 & 3 & 4 & 5}{-+++} &
        \stab{i}{3 \\ 1 & 2 & 4 & 5}{+-++} &
        \stab{j}{4 \\ 1 & 2 & 3 & 5}{++-+} &
        \stab{k}{5 \\ 1 & 2 & 3 & 4}{+++-}
      \end{array} \\[2\vsp]
      \psset{linewidth=.1ex}
      \everypsbox{\scriptstyle}
      \ncline[nodesep=3pt,linecolor=red]    {h}{i} \naput{\red 2}
      \ncline[nodesep=3pt,linecolor=blue]   {i}{j} \naput{\blue 3}
      \ncline[nodesep=3pt,linecolor=green]  {j}{k} \naput{\green 4}
        \G_{(3,2)} &
      \begin{array}{\cci \cci \cci \cci \cci}
        \stab{a}{3 & 4 \\ 1 & 2 & 5}{+-++} &
        \stab{b}{2 & 4 \\ 1 & 3 & 5}{-+-+} &
        \stab{c}{2 & 5 \\ 1 & 3 & 4}{-++-} &
        \stab{d}{3 & 5 \\ 1 & 2 & 4}{+-+-} &
        \stab{e}{4 & 5 \\ 1 & 2 & 3}{++-+}
      \end{array} \\[2\vsp]
      \psset{linewidth=.1ex}
      \everypsbox{\scriptstyle}
      \ncline[nodesep=3pt,offset=2pt,linecolor=red]  {a}{b} \naput{\red 2}
      \ncline[nodesep=3pt,offset=2pt,linecolor=blue] {b}{a} \naput{\blue 3}
      \ncline[nodesep=3pt,linecolor=green]           {b}{c} \naput{\green 4}
      \ncline[nodesep=3pt,linecolor=red]             {c}{d} \naput{\red 2}
      \ncline[nodesep=3pt,offset=2pt,linecolor=blue] {d}{e} \naput{\blue 3}
      \ncline[nodesep=3pt,offset=2pt,linecolor=green]{e}{d} \naput{\green 4}
        \G_{(3,1,1)} &
      \begin{array}{\cci \cci \cci \cci \cci}
        & & \stab{w}{4 \\ 3 \\ 1 & 2 & 5}{+--+} & & \\[-1\vsp]
        \stab{u}{3 \\ 2 \\ 1 & 4 & 5}{--++} &
        \stab{v}{4 \\ 2 \\ 1 & 3 & 5}{-+-+} & &
        \stab{y}{5 \\ 3 \\ 1 & 2 & 4}{+-+-} &
        \stab{z}{5 \\ 4 \\ 1 & 2 & 3}{++--} \\[-1\vsp]
        & & \stab{x}{5 \\ 2 \\ 1 & 3 & 4}{-++-} & &
      \end{array}
      \psset{linewidth=.1ex}
      \everypsbox{\scriptstyle}
      \ncline[nodesep=3pt,linecolor=blue]   {u}{v}  \naput{\blue 3}
      \ncline[nodesep=3pt,linecolor=red]    {v}{w} \naput{\red 2}
      \ncline[nodesep=3pt,linecolor=green]  {v}{x}  \nbput{\green 4}
      \ncline[nodesep=3pt,linecolor=green]  {w}{y} \naput{\green 4}
      \ncline[nodesep=3pt,linecolor=red]    {x}{y}  \nbput{\red 2}
      \ncline[nodesep=3pt,linecolor=blue]   {y}{z}  \naput{\blue 3}
    \end{array}
  \end{displaymath}
    \caption{\label{fig:G5}The standard dual equivalence graphs
      $\G_{(4,1)}, \G_{(3,2)}$ and $\G_{(3,1,1)}$.}
  \end{center}
\end{figure}

In contrast with $\X_{\lambda}^{n}$, the bound on the edge colors
(equivalently, on the entries of the tableaux) for $\G_{\lambda}$ is
implicit. More precisely, if $\lambda$ is a partition of $m$, then
$\G_{\lambda}$ has degree $m$, meaning that the largest edge color is
$m\!-\!1$.

It follows from results in \cite{Haiman1992} that each $\G_{\lambda}$
is connected. It is easy to see from the construction that
$\G_{\lambda'}$ can be obtained from $\G_{\lambda}$ by conjugating
each vertex and multiplying each entry of the signature by
$-1$. Several other nice properties of these graphs are given
\cite{Assaf2007-2}, including the fact that the $\G_{\lambda}$ are
pairwise non-isomorphic and have no nontrivial automorphisms.

Motivated by the local characterization of crystal graphs
$\X_{\lambda}^{n}$, we next present a local characterization of the
graphs $\G_{\lambda}$. We will abuse notation by simultaneously
referring to $D_i$ as the collection of $i$-colored edges as well as
an involution on vertices which have an $i$-edge, i.e. we write
$D_i(v)=w$ whenever $\{v,w\} \in D_i$. This is justified by the first
axiom.

\begin{definition}[\cite{Assaf2007-2}]
  A signed, colored graph $\G = (V,\sigma,D)$ is a {\em dual
    equivalence graph} if the following hold:
  \begin{itemize}
    
  \item[(ax1)] For $v \in V$ and $i>1$, $\sigma(v)_{\imo} =
    -\sigma(v)_{i}$ if and only if there exists $w \in V$ such that
    $\{v,w\} \in D_{i}$.  Moreover, $w$ is unique when it exists.

  \item[(ax2)] Whenever $\{v,w\} \in D_{i}$,
      \begin{displaymath}
        \begin{array}{rcrl}
          \sigma(v)_j & = & -\sigma(w)_j & \; \mbox{for} \; j=\imo,i ;\\
          \sigma(v)_h & = &  \sigma(w)_h & \; \mbox{for} \; h < \imt \;\;
          \mbox{and} \;\; h > \ipo . \\
        \end{array}
      \end{displaymath}
      
  \item[(ax3)] Whenever $\{v,w\} \in D_{i}$,
    \begin{displaymath}
      \begin{array}{ll}
        \mbox{if} \; \sigma(v)_{\imt} = -\sigma(w)_{\imt}, & 
        \mbox{then} \; \sigma(v)_{\imt} = -\sigma(v)_{\imo}; \\
        \mbox{if} \; \sigma(v)_{\ipo} = -\sigma(w)_{\ipo}, &
        \mbox{then} \; \sigma(v)_{\ipo} = -\sigma(v)_{i}.
      \end{array}
    \end{displaymath}
    
  \item[(ax4)] For $i>3$, every non-trivial connected component of
    $(V,\sigma,D_{\imt} \cup D_{\imo} \cup D_{i})$ is depicted below
    and there is a path containing at most one $D_i$ edge between any
    two vertices of $(V,\sigma,D_2 \cup \cdots \cup D_i)$.
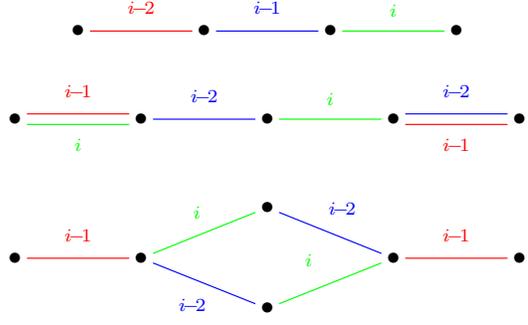
    \begin{displaymath}
      \begin{array}{c}
        \begin{array}{\cciii \cciii \cciii \cciii}
          \rn{h}{\B} & \rn{i}{\B} & \rn{j}{\B} & \rn{k}{\B}
        \end{array} \\[1.5\vsp]
        \begin{array}{\cciii \cciii \cciii \cciii \cciii}
          \rn{a}{\B} & \rn{b}{\B} & \rn{c}{\B} & \rn{d}{\B} & \rn{e}{\B}
        \end{array} \\[1.5\vsp]
        \begin{array}{\cciii \cciii \cciii \cciii \cciii}
          & & \rn{w}{\B} & & \\[.5\vsp]
          \rn{u}{\B} & \rn{v}{\B} & & \rn{y}{\B} & \rn{z}{\B} \\[.5\vsp]
          & & \rn{x}{\B} & &
        \end{array}
        \psset{nodesep=2pt,linewidth=.1ex}
        \everypsbox{\scriptstyle}
        \ncline[linecolor=red]    {h}{i} \naput{\red  \imt}
        \ncline[linecolor=blue]   {i}{j} \naput{\blue \imo}
        \ncline[linecolor=green]  {j}{k} \naput{\green i}
        \ncline[linecolor=red,offset=2pt]  {a}{b} \naput{\red \imo}
        \ncline[linecolor=green,offset=2pt]{b}{a} \naput{\green i }
        \ncline[linecolor=blue]            {b}{c} \naput{\blue \imt}
        \ncline[linecolor=green]           {c}{d} \naput{\green i }
        \ncline[linecolor=blue,offset=2pt] {d}{e} \naput{\blue \imt}
        \ncline[linecolor=red,offset=2pt]  {e}{d} \naput{\red \imo}
        \ncline[linecolor=red]   {u}{v}\naput{\red \imo}
        \ncline[linecolor=green] {v}{w}\naput{\green i }
        \ncline[linecolor=blue]  {v}{x}\nbput{\blue \imt}
        \ncline[linecolor=blue]  {w}{y}\naput{\blue \imt}
        \ncline[linecolor=green] {x}{y}\naput{\green \!\!\!\!\!\!i}
        \ncline[linecolor=red]   {y}{z}\naput{\red \imo}
      \end{array}
    \end{displaymath}

   \item[(ax5)] Whenever $|i-j| \geq 3$, $\{v,u\} \in D_i$ and
     $\{u,w\} \in D_j$, there exists $x \in V$ such that $\{v,x\} \in
     D_j$ and $\{x,w\} \in D_i$.
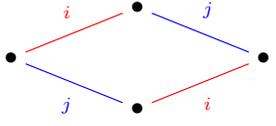
     \begin{displaymath}
       \begin{array}{\cciii\cciii\cciii}
                   &\rn{x}{\B}&           \\[.5\vsp]
         \rn{w}{\B}&          &\rn{y}{\B} \\[.5\vsp]
                   &\rn{v}{\B}& 
       \end{array}
       \psset{nodesep=3pt,linewidth=.1ex}
       \everypsbox{\scriptstyle}
       \ncline[linecolor=red]{w}{x}  \naput{\red i}
       \ncline[linecolor=blue]{w}{v} \nbput{\blue j}
       \ncline[linecolor=blue]{x}{y} \naput{\blue j}
       \ncline[linecolor=red]{v}{y}  \nbput{\red i}
     \end{displaymath}
  \end{itemize}
\label{defn:deg}
\end{definition}

When $n = 4$, we stipulate the following in place of axiom $4$: For
$i>2$, every non-trivial connected component of the sub-graph
$(V,\sigma, D_{\imo} \cup D_{i})$ is one of the following:
    
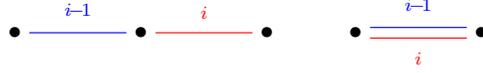
\begin{displaymath}
      \begin{array}{\cciii\cciii\ccii\cciii\cciii}
        \\[.3\vsp]
        \rn{a}{\B} & \rn{b}{\B} & \rn{c}{\B} & \rn{d}{\B} & \rn{e}{\B}
        \\[.6\vsp] 
      \end{array}
      \psset{nodesep=3pt,linewidth=.1ex}
      \everypsbox{\scriptstyle}
      \ncline[linecolor=blue]           {a}{b} \naput{\blue \imo}
      \ncline[linecolor=red]           {b}{c} \naput{\red i}
      \ncline[linecolor=blue,offset=2pt]{d}{e} \naput{\blue \imo}
      \ncline[linecolor=red,offset=2pt]{e}{d} \naput{\red i}
    \end{displaymath}

Next we have the analog of \refthm{structure-crystal} for dual
equivalence graphs. The proof of \refthm{isomorphic} is by a purely
combinatorial, graph-theoretic argument.

\begin{theorem}[\cite{Assaf2007-2}]
  Each $\G_{\lambda}$ is a dual equivalence graph, and every connected
  component of a dual equivalence graph is isomorphic to
  $\G_{\lambda}$ for a unique partition $\lambda$.
\label{thm:isomorphic}
\end{theorem}

%
\section{From crystal graphs to dual equivalence graphs}
%
\label{sec:correspondence}

Since standard tableaux are special cases of semi-standard tableaux,
and since the construction of $\G_{\lambda}$ was motivated by
$\X_{\lambda}^{n}$, we would like to explore the connection between
regular graphs and dual equivalence graphs. Since semi-standard
tableaux index a basis for an irreducible representation of $\SL_n$ so
that the content of the tableau corresponds to the weight of the basis
vector, those basis vectors corresponding to standard tableaux index
the $0$-weight space of the representations of dimension
$n$. Therefore we may think of dual equivalence graphs as living in
the $0$-weight space of crystal graphs.

In general, the Weyl group acts naturally on the $0$-weight space of a
representation. In the case of $\SL_n$, the Weyl group is isomorphic
to $\Sn$. Since standard tableaux index the basis for an irreducible
representation of $\Sn$, we would like to use dual equivalence graphs
to combinatorialize the study of representations of $\Sn$ in much the
same way that crystal graphs combinatorialize the study of
representations of $\SL_n$. This is the motivation behind the
following correspondence. In the interest of space, the longer proofs
have been condensed by omitting the more mechanical details.

Combinatorially, the $0$-weight space of a regular graph $\X$ consists
of all vertices which lie at the center of each $i$-string;
i.e. $\varepsilon(x,i) = -\delta(x,i)$ for each $i$. Given a directed,
colored graph $\X$, define a signed, colored graph $\G(\X) =
(V,\sigma,D)$ by
\begin{eqnarray}
  V & = & \left\{ x \in \X \; | \varepsilon(x,i) = -\delta(x,i) = 0
    \;\mbox{or}\; 1 \; \;\forall\; i, \right\}, 
  \label{eqn:vertices}\\
  \sigma(x)_i & = & \left\{ \begin{array}{rl}
      +1 & \mbox{if}\; \varepsilon(x,i) = 1 \\
      -1 & \mbox{if}\; \varepsilon(x,i) = 0
    \end{array} \right. ,
  \label{eqn:sigma}\\
  D_{i}(x) & = & \left\{ \begin{array}{rl}
      F_{\imo}F_iE_{\imo}E_i x & \mbox{if}\; \varepsilon(x,i) = 1
      \;\mbox{and}\; \varepsilon(x,\imo) = 0 \\
      F_iF_{\imo}E_iE_{\imo} x & \mbox{if}\; \varepsilon(x,i) = 0
      \;\mbox{and}\; \varepsilon(x,\imo) = 1
    \end{array} \right. .
  \label{eqn:edges}
\end{eqnarray}
An $i$-edge from $x$, or equivalently $D_i(x)$, is not defined if
neither condition of \refeq{edges} holds.

\begin{figure}[ht]
  \begin{center}
    \begin{displaymath}
      \begin{array}{\cci\cci\cci\cci\cci}
        & & \rn{a11}{\B} & & \\[1.5\vsp]
        & \rn{b10}{\B} & & \rn{b12}{\B} & \\[1.5\vsp]
        & & \begin{array}{c} 
          \rn{c11l}{\B} \\[-.3\vsp] _{+-} 
        \end{array} \hspace{4ex} \begin{array}{c} 
          \rn{c11r}{\B} \\[-.3\vsp] _{-+} 
        \end{array} & & \\[.7\vsp] 
        \rn{d9}{\B} & & & & \rn{d13}{\B} \\[.7\vsp]
        & & \rn{e11}{\B} & & 
      \end{array}
      \psset{linewidth=.1ex,nodesep=3pt}
      \everypsbox{\scriptstyle}
      \ncline[linewidth=.3ex,linestyle=dashed] {c11l}{c11r} \nbput{D_i}
      \ncline[linecolor=blue]{->} {a11}{b10}  
      \ncline[linewidth=.3ex,linecolor=red]{->} {a11}{b12}   
      \ncline[linewidth=.3ex,linecolor=red]{->} {b10}{c11r}  
      \ncline[linecolor=blue]{->} {b12}{c11l} 
      \ncline[linecolor=blue]{->} {c11l}{d9}  
      \ncline[linewidth=.3ex,linecolor=red]{->} {c11r}{d13}  
      \ncline[linewidth=.3ex,linecolor=red]{->} {d9}{e11}    
      \ncline[linecolor=blue]{->} {d13}{e11}  
    \end{displaymath}      
    \caption{\label{fig:edges}Constructing $\G(\X)$ from $\X$; here we
    have $F_{\imo} {\blue \swarrow}$, $F_{i} {\red \searrow}$.}
  \end{center}
\end{figure}

\begin{proposition}
  If $\X$ is a regular graph, then $\G(\X)$ is well-defined.
\end{proposition}

\begin{proof}
  It suffices to show that $D_i$ is well-defined. Note that if
  $\delta(x,i)=-1$, then $E_i x$ is nonzero. If $\delta(x,\imo) = 0$,
  then by axiom P3 we must have $\delta(E_i x, \imo) + \varepsilon(E_i
  x, \imo) = -1$. The constraints of axiom P4 force $\varepsilon(E_i
  x, \imo) = 0$, and so we conclude that $\delta(E_i x, \imo) =
  -1$. Therefore $E_{\imo} E_i x$ is nonzero. Continuing thus, with
  liberal use of axioms P3 and P4, we eventually conclude that
  $F_{\imo}F_iE_{\imo}E_i x$ is nonzero. The analogous argument holds
  with the reverse assumptions for $\delta(x,\imo), \delta(x,i)$.
\end{proof}

Below we give two separate proofs that if $\X$ is a regular graph with
$0$-weight space given by $V$, then $\G(\X)$ is a dual equivalence
graph. The first proof utilizes the graph structures on tableaux, and
the second uses the local characterizations. Note that the condition
that $V$ give the $0$-weight space is highly restrictive. 

In order for $\X_{\lambda}^{n}$ to have a $0$-weight space, $\lambda$
must be a partition of $kn$ into at most $n$ parts. We first treat the
case for $\X_{\lambda}^{n}$ when $k=1$, i.e. $\lambda$ is a partition
of $n$.

\begin{theorem}
  For every partition $\lambda$ of $n$, we have $\G(\X_{\lambda}^{n})
  = \G_{\lambda}$.
\label{thm:tableaux-pf}
\end{theorem}

\begin{proof}
  In this case we may identify $V$ with $\SYT(\lambda)$. Let $w_T$ be
  the reading word for a tableau $T$. For $T \in \SYT(\lambda)$, note
  that $\varepsilon(T,i) = 0$ if and only if $\ipo$ occurs to the left
  of $i$ in $w_T$. Therefore equations \ref{eqn:sigma-T} and
  \ref{eqn:sigma} correspond. In particular, $T$ will have an $i$-edge
  in $\G_{\lambda}$ precisely when \refeq{edges} is defined.

  Suppose $D_i(T)$ is defined, and so $T$ admits an elementary dual
  equivalence for $\imo,i,\ipo$. By symmetry, we may assume $i$ lies
  to the left of $i\pm 1$ in $w_T$. If $\varepsilon(T,i)=1$, then
  ignoring all letters other than $\imo,i,\ipo$ and regarding
  $\we_j,\wf_j$ ($j=\imo,i$) as operators on reading words, we have
  \begin{displaymath}
    \begin{array}{r c r ccccccl}
      \wf_{\imo}\wf_{i}\we_{\imo}\we_{i} w_T
      & = & \wf_{\imo}\wf_{i}\we_{\imo}\we_{i}  & \cdots & i & \cdots &
      \makebox[0pt]{$\imo$} & \cdots & \makebox[0pt]{$\ipo$} & \cdots \\
      & = & \wf_{\imo}\wf_{i}\we_{\imo}  & \cdots & i & \cdots &
      \makebox[0pt]{$\imo$} & \cdots & i & \cdots \\ 
      & = & \wf_{\imo}\wf_{i}  & \cdots & i & \cdots & \makebox[0pt]{$\imo$} &
      \cdots & \makebox[0pt]{$\imo$} & \cdots \\
      & = & \wf_{\imo}  & \cdots & \makebox[0pt]{$\ipo$} & \cdots &
      \makebox[0pt]{$\imo$} & \cdots & \makebox[0pt]{$\imo$} & \cdots \\ 
      & = & & \cdots & \makebox[0pt]{$\ipo$} & \cdots &
      \makebox[0pt]{$\imo$} & \cdots & i & \cdots . 
    \end{array}
  \end{displaymath}
  Therefore \refeq{edges} corresponds to the elementary dual
  equivalences. The case when $\varepsilon(T,\imo)=1$ is similar.
\end{proof}

Going back to the representation theory of $\SL_n$ and $\Sn$,
\refthm{tableaux-pf} reflects the fact that for $\lambda$ a partition
of $n$, the action of $\Sn$ on the $0$-weight space of the irreducible
representation of $\SL_n$ corresponding to $\lambda$ is given by the
Specht module $S^{\lambda}$.

The following result is a corollary to \refthm{tableaux-pf}, assuming
Theorems \ref{thm:structure-crystal} and \ref{thm:isomorphic}. Here we
present a self-contained proof based solely on the local
characterizations of crystal graphs and dual equivalence graphs
without appealing to Theorems \ref{thm:structure-crystal} and
\ref{thm:isomorphic} and without reference to tableaux.

\begin{theorem}
  If $\X$ is a regular with $0$-weight space given by $V$ (see
  \refeq{vertices}), then $\G(\X)$ is a dual equivalence graph.
\label{thm:axiom-pf}
\end{theorem}

\begin{proof}
  Comparing \refeq{sigma} with the conditions in \refeq{edges} 
  dictating when $D_i$ is defined shows that $\sigma(x)_{\imo} =
  -\sigma(x)_i$ if and only if $D_i(x)$ is defined. \reffig{edges}
  makes evident the fact that $\sigma(x)_j = -\sigma(D_i(x))_j$ for $j
  = \imo,i$, though this can also be seen directly from the proof that
  $D_i$ is well-defined. Therefore we have established dual
  equivalence axioms 1 and 2.

  Axiom 3 can be demonstrated by showing that if $\varepsilon(x,\imt)
  \neq \varepsilon(D_i(x),\imt)$, then $D_i(x) = D_{\imo}(x)$. To do
  this, we show that the crystal graph operators going from $x$ to
  $D_i(x)$ exactly mirror the situation in $\X_{(2,2)}^{4}$ for the
  two standard tableaux of shape $(2,2)$; see \reffig{X22}.

  By crystal axioms P3, P4 and P5, whenever $|a-b| \geq 2$, we have
  $E_a E_b = E_b E_a$, $F_a E_b = E_b F_a$ and $F_a F_b = F_b
  F_a$. Dual equivalence axiom 5 is equivalent to showing that if
  $D_i$ and $D_j$ are both defined at $x$ for some $|i-j| \geq 3$,
  then $D_i D_j (x) = D_j D_i (x)$. Therefore the identities above for
  crystal operators establish dual equivalence axiom 5.

  \begin{figure}[ht]
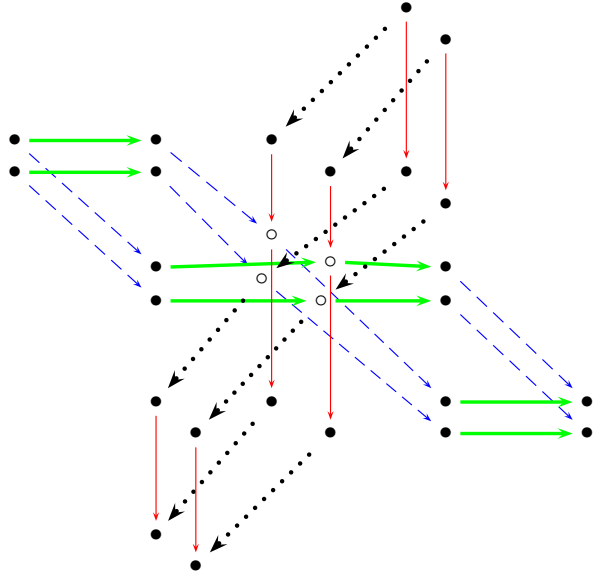

    \begin{center}
      \begin{displaymath}
        \begin{array}{\ccii cccc cccc ccc\ccii c}
          & & & & & & & & & \rn{a10}{\B} & & & & \\
          & & & & & & & & & & \rn{b11}{\B} & & & \\
          \\[\vsp]
          \rn{d01}{\B} & & & \rn{d04}{\B} & & & \rn{d07}{\B} & & & & &
          & & \\
          \rn{e01}{\B} & & & \rn{e04}{\B} & & & & \rn{e08}{\B} & &
          \rn{e10}{\B} & & & & \\          
          & & & & & & & & & & \rn{f11}{\B} & & & \\
          & & & & & & \rn{g07}{\circ} & & & & & & & \\
          & & & \rn{h04}{\B} & & & \raisebox{-1ex}{$\rn{h07}{\circ}$}
          \hspace{1.7ex} & \raisebox{.5ex}{$\rn{h08}{\circ}$} 
          & & & \rn{h11}{\B} & & & \\ 
          & & & \rn{i04}{\B} & & & & \rn{i08}{\circ} \hspace{1.7ex} & & &
          \rn{i11}{\B} & & & \\
          \\[\vsp]
          & & & \rn{k04}{\B} & & & \rn{k07}{\B} & & & & \rn{k11}{\B} &
          & & \rn{k14}{\B} \\
          & & & & \rn{l05}{\B} & & & \rn{l08}{\B} & & & \rn{l11}{\B} &
          & & \rn{l14}{\B} \\
          \\[\vsp]
          & & & \rn{n04}{\B} & & & & & & & & & & \\
          & & & & \rn{o05}{\B} & & & & & & & & &
        \end{array}
        \psset{linewidth=.1ex,nodesep=3pt}
        \everypsbox{\scriptstyle}
        \ncline[linecolor=blue,linestyle=dashed]{->}{d01}{h04}
        \ncline[linecolor=blue,linestyle=dashed]{->}{e01}{i04}
        \ncline[linecolor=blue,linestyle=dashed]{->}{d04}{g07}
        \ncline[linecolor=blue,linestyle=dashed]{->}{e04}{h07}
        \ncline[linecolor=blue,linestyle=dashed]{->}{g07}{k11}
        \ncline[linecolor=blue,linestyle=dashed]{->}{h07}{l11}
        \ncline[linecolor=blue,linestyle=dashed]{->}{h11}{k14}
        \ncline[linecolor=blue,linestyle=dashed]{->}{i11}{l14}
        \ncline[linecolor=green,linewidth=.3ex]{->}{d01}{d04}
        \ncline[linecolor=green,linewidth=.3ex]{->}{e01}{e04}
        \ncline[linecolor=green,linewidth=.3ex]{->}{h04}{h08}
        \ncline[linecolor=green,linewidth=.3ex]{->}{i04}{i08}
        \ncline[linecolor=green,linewidth=.3ex]{->}{h08}{h11}
        \ncline[linecolor=green,linewidth=.3ex]{->}{i08}{i11}
        \ncline[linecolor=green,linewidth=.3ex]{->}{k11}{k14}
        \ncline[linecolor=green,linewidth=.3ex]{->}{l11}{l14}
        \ncline[linecolor=black,linewidth=.4ex,linestyle=dotted]{->}{a10}{d07}
        \ncline[linecolor=black,linewidth=.4ex,linestyle=dotted]{->}{b11}{e08}
        \ncline[linecolor=black,linewidth=.4ex,linestyle=dotted]{->}{e10}{h07}
        \ncline[linecolor=black,linewidth=.4ex,linestyle=dotted]{->}{f11}{i08}
        \ncline[linecolor=black,linewidth=.4ex,linestyle=dotted]{->}{h07}{k04}
        \ncline[linecolor=black,linewidth=.4ex,linestyle=dotted]{->}{i08}{l05}
        \ncline[linecolor=black,linewidth=.4ex,linestyle=dotted]{->}{k07}{n04}
        \ncline[linecolor=black,linewidth=.4ex,linestyle=dotted]{->}{l08}{o05}
        \ncline[linecolor=red]{->}{a10}{e10}
        \ncline[linecolor=red]{->}{b11}{f11}
        \ncline[linecolor=red]{->}{d07}{g07}
        \ncline[linecolor=red]{->}{e08}{h08}
        \ncline[linecolor=red]{->}{g07}{k07}
        \ncline[linecolor=red]{->}{h08}{l08}
        \ncline[linecolor=red]{->}{k04}{n04}
        \ncline[linecolor=red]{->}{l05}{o05}
      \end{displaymath}      
      \caption{\label{fig:flower} An illustration of the commutativity
        which establishes dual equivalence axiom 4; the $0$-weight
        vertices are $\circ$'s and crystal operators are $F_{\imh}
        {\blue \searrow}$, $F_{\imt} {\green \rightarrow}$, $F_{\imo}
        {\swarrow}$, $F_{i} {\red \downarrow}$.}
    \end{center}
  \end{figure}

  The most difficult step is to establish axiom 4. For this, note that
  it is sufficient to show that if $D_{\imt}, D_{\imo}, D_{i}$ are all
  defined at $x$ and are all distinct (i.e. $x$ does not have a double
  edge), then $D_{\imt} D_{i}(x) = D_{i} D_{\imt}(x)$. By symmetry we
  may assume $\varepsilon(x,\imh) = 1 = \varepsilon(x,\imo)$ and
  $\varepsilon(x,\imt) = 0 = \varepsilon(x,i)$. Then the first part of
  axiom 4 is equivalent to establishing the following identity,
  illustrated in \reffig{flower}.
  \begin{displaymath}
    (F_{i} F_{\imo} E_{i} E_{\imo}) (F_{\imt} F_{\imh} E_{\imt}
    E_{\imh}) x \; = \; (F_{\imt} F_{\imh} E_{\imt} E_{\imh}) (F_{i}
    F_{\imo} E_{i} E_{\imo}) x
  \end{displaymath}
  The key to establishing this relation is to use the commutativity
  relations for crystal operators used to prove dual equivalence axiom
  5, as well as the additional relation that both $E_{\imo}$ and
  $F_{\imo}$ commute with the sequence $F_{\imt}F_{\imh}E_{\imt}$.

  For the second part of axiom 4, we first show that if $\X$ is
  connected, then so is $\G(\X)$. In particular, we show that any two
  vertices of $V$ may be connected by a path consisting of sequences
  as in \refeq{edges}. Furthermore, for any two vertices we may choose
  a path containing at most one sequence using the maximal edge
  color. Now axiom 4 follows.
\end{proof}

\begin{corollary}
  Let $\lambda = (\mu_1 + k, \mu_2 + k, \ldots, \mu_n + k)$ for some
  partition $\mu$ of $n$ and some non-negative integer $k$. Then
  $\G(\X_{\lambda}^{n}) = \G_{\mu}$.
\label{cor:addcol}
\end{corollary}

The partition $\lambda$ is the result of adding $k$ columns of height
$n$ to $\mu$. The graphs $\X_{\lambda}^{n}$ for such $\lambda$ give
all of the regular graphs which have $0$-weight spaces given by
$V$. Such a graph corresponds to the irreducible representation of
$\SL_n$ indexed by $\mu$, since $\SL_n$ cannot discriminate between
the irreducible representations indexed by $\lambda$ and
$\mu$. Therefore \refcor{addcol} shows that the corresponding dual
equivalence graph decomposes (albeit trivially) the action of the Weyl
group on the $0$-weight space of the representation.

\begin{remark}
  By introducing weights to the regular graphs, we can recover the
  {\em dimension} of a regular graph, which in the context of
  \refthm{structure-crystal} gives the size of the indexing
  partition. Letting $k$ denote the ratio of the dimension of $\X$ to
  the degree of $\X$, modify the definition of $\sigma$ in \refeq{sigma}
  to be
  \begin{equation}
    \sigma(x)_i = \left( -1 \right)^{\varepsilon(x,i) + k}.
    \label{eqn:sigmak}
  \end{equation}
  Since multiplying the signatures of a dual equivalence graph
  entry-wise by $-1$ has the effect of conjugating the indexing
  partition, under this definition \refcor{addcol} must read
  \begin{equation}
    \G(\X_{\lambda}^{n}) = \left\{ \begin{array}{ll}
        \G_{\mu} & \mbox{if $k$ is odd}, \\
        \G_{\mu'} & \mbox{if $k$ is even}.
      \end{array} \right.
  \end{equation}
  Representation theoretically, this corresponds to tensoring with the
  determinant representation in $\GL_n$, which has the effect of
  conjugating the shape. So indeed, the dual equivalence graph on the
  $0$-weight space of $\X_{\lambda}^{n}$ decomposes the action of
  $\Sn$ on the $0$-weight space of the corresponding representation of
  $\GL_n$ with this minor modification of signatures.
\end{remark}

%
\section{Extensions}
%
\label{sec:extensions}

The ultimate goal is to extend the construction of $\G(\X)$ for any
regular graph $\X$ which has a $0$-weight space. The characterization
of $V$ in \refeq{vertices} would become
\begin{equation}
  V = \left\{ x \in \X \; | \; \varepsilon(x,i) = -\delta(x,i)
    \;\forall\; i \right\}.
\label{eqn:vertices2}
\end{equation}
In many cases, the $0$-weight graphs obtained from regular graphs in a
similar manner can be viewed as D graphs, which are generalizations of
dual equivalence graphs; see \cite{Assaf2007-2}. Each D graph can be
transformed into a disjoint union of dual equivalence graphs, so the
hope is that each $0$-weight graph may be realized as a D graph which
decomposes into dual equivalence graphs in the correct way. If
successful, this would provide a combinatorial way to study the
representations of Weyl groups on $0$-weight spaces, at least in type
A.

A further direction is to use the correspondence between type A
crystal graphs and dual equivalence graphs to define dual equivalence
graphs in other types. For the simply-laced types at least, this can
be done by the same constructions using Stembridge's axioms. Using the
combinatorial description of crystal bases in terms of generalized
Young tableaux \cite{KaNa1994,Littelmann1995}, we may define the
standard graphs which are the analogs of $\G_{\lambda}$, and then use
the local characterization of crystal graphs to obtain a local
characterization for dual equivalence graphs.

%
%

\bibliographystyle{amsalpha} 
\bibliography{../references}

\end{document}